\documentclass[12pt,reqno]{amsart}

\usepackage{amsmath,amsthm}
\usepackage{tikz}
\usepackage{amssymb}
\usepackage[english]{babel}

\usepackage{subcaption}
\usepackage{times}
\usepackage[margin=1.5in]{geometry}
\usepackage{mathrsfs}
\usepackage{hyperref}
\usepackage{xcolor}

\newcommand{{\rr}}{\mathcal{R}}

\newtheorem{thm}{Theorem}[section]
\newtheorem{corr}[thm]{Corollary}
\newtheorem{lem}[thm]{Lemma}

\newtheorem{ques}[thm]{Question}

\newcommand{\pf}{\noindent\textbf{Proof.}$\ $}
\newcommand{\zb}{$\hfill\Box$}
\allowdisplaybreaks

\makeatletter\@addtoreset{equation}{section} \makeatother

\begin{document}

	\title [The density of primes in the eigensurface of ${\bf S}_3$]{The density of primes in the eigensurface of the symmetric group ${\bf S}_3$}

	\author[L. Geng]{Liang Geng}
	\footnotemark[1]
	\address{Liang Geng: School of Mathematics, Southeast University,
		Nanjing, Jiangsu 211189, China.} \email{lgeng@seu.edu.cn}

	\author[W. He]{Wei He$^*$}
	\thanks{$^*$ corresponding author}
	\address{Wei He: School of Mathematics, Southeast University,
		Nanjing, Jiangsu 211189, China.} \email{hewei@seu.edu.cn}
	
	\author[R. Yang]{Rongwei Yang}
	\address{Rongwei Yang: Department of Mathematics and Statistics, University at Albany, the State University of New York,
		Albany, NY 12222, U.S.A.} \email{ryang@albany.edu}

	\maketitle
	
\begin{abstract}
The Prime Number Theorem asserts that the density of primes less than or equal to $N$ is asymptotically equal to $1/\log N$. The density of prime triples in coprime triples in $\mathbb{Z}^3_+$ is determined to be $3\zeta (3)/\log N$, where $\zeta$ is the Riemann zeta function. In this paper, we prove that the density of prime triples in coprime triples in the surface ${\mathcal S}=\{z_0^{2} - z_1^{2} + z_2^{2} - z_0z_2=0\}$ is greater than $3\zeta (3)/\log N$, meaning that $\mathcal{S}$ meets primes more frequently. This surface is the eigensurface of the symmetric group ${\bf S}_3$ with respect to an irreducible representation.\\

		
\noindent Key words and phrases: eigensurface, density of primes, prime triples, symmetric group, dihedral group. 
		
\noindent Mathematics Subject Classification (2020): 11A41; 11D45; 47A13; 20C30

\end{abstract}

\section{Introduction}

	Let $\mathbb{Z}_+$ be the set of nonnegative integers. In this paper, a triple $z=(z_0,z_1,z_2)\in \mathbb{Z}^3_+$ is said to be {\em coprime} (or {\em relatively prime}) if the greatest common divisor of $z_0,z_1,z_2$, denoted by $\gcd (z)$, is equal to $1$. A triple $z$ is called {\em prime} if it is coprime, and at least one component of $z$ is prime, for example $(3, 8, 8)$ and $(5, 0, 6)$. This paper studies the density of prime triples in coprime triples in algebraic surfaces in $\mathbb{C}^3$. Clearly, the majority of such surfaces don't contain any non-zero integer triples, for instance the surface $\{ z_0^2=\sqrt{2}z_1z_2\}$ or $\{z_0^3+z_1^3=z_2^3\}$, not mentioning prime triples. Thus, to make the study meaningful, it is important to pick a ``right" class of surfaces. This paper considers the surface ${\mathcal S}=\{z_0^{2} - z_1^{2} + z_2^{2} - z_0z_2=0\}$ which is the {\em projective spectrum} (or {\em eigensurface}) associated with an irreducible representation of the symmetric group ${\bf S}_3$.  Although this surface is specific, the method used here can be generalized to work for some other algebraic surfaces.	
	 In the sequel, we review relevant concepts and state the main theorem of this paper.
	
	
The notion of projective spectrum was introduced by the third author \cite{Ya1}. Let $\mathcal{B}$ be a unital Banach algebra, and $\mathbb{A}= (A_0, A_1, \dots, A_n)$ be a tuple of elements in $\mathcal{B}$. For $z= (z_0, \dots, z_n) \in \mathbb{C}^{n+1}$, we set $\mathbb{A}(z) = z_0A_0 + \dots + z_nA_n$. The {\em projective spectrum} of $\mathbb{A}$ is defined as \[P(\mathbb{A})= \{z \in \mathbb{C}^{n+1} \ |\ \mathbb{A}(z) \text  { is not invertible in $\mathcal{B}$}\}.\] Evidently, when $\mathcal{B}$ is a matrix algebra, the projective spectrum $ P(\mathbb{A})$ is the algebraic variety $\{\det \mathbb{A}(z)= 0\}$. In particular, the polynomial $Q_{\mathbb{A}}(z):=\det (z_0I+z_1A_1+\cdots +z_n A_n)$ is called the ({\em joint}) {\em characteristic polynomial} of the matrices, and its zero variety $Z(Q_{\mathbb{A}})$ is called the corresponding {\em eigensurface}. In the case where $A_1, \dots, A_n$ are representations of a generating set of a group $G$, it is known that the projective spectrum $P({\mathbb A})$ captures some deep intrinsic properties of $G$ and the representation. For example, for certain self-similar groups, such as the infinite dihedral group, the lamplighter group, and the Grigorchuk group of intermediate growth, the projective spectra of ${\mathbb A}$ are shown to be closely related to the Julia sets of spectral dynamic systems for the groups \cite{GY, GrY, ZYL}. In a different direction, the characteristic polynomials for finite dimensional Lie algebras have also been actively studied \cite{AY, GLW, HuY}. For more details about these developments, we refer the reader to monograph \cite{Ya2} and the references therein.

This paper aims to initiate an investigation on the number theoretic properties of eigensurfaces associated with finite groups. The focus here is the symmetric group ${\bf S}_3$, which is isomorphic to the dihedral group $D_3=\langle a, t\mid a^2=t^2=(at)^3=1\rangle$. An irreducible unitary representation of $D_3$ is defined by the map 
\[	\rho(a)=
	\begin{pmatrix}
		0 & e^{2\pi i/3} \\
		e^{-2\pi i/3} & 0 
	\end{pmatrix},
	\quad \rho(t) =
	\begin{pmatrix}
		0 & 1 \\
		1 & 0 
	\end{pmatrix}.\]
Then the characteristic polynomial for $(\rho(a), \rho(at))$ is
\begin{align*}
	Q(z)=\det (z_0I+z_1\rho(a)+z_2\rho(at))=z_0^{2} - z_1^{2} + z_2^{2} - z_0z_2,
\end{align*}
and the corresponding eigensurface is 
\begin{align}\label{surfaceequ}
	{\mathcal S}=\{z_0^{2} - z_1^{2} + z_2^{2} - z_0z_2=0\}.
\end{align}

It is worth noting that not every finite group has an eigensurface containing prime tuples, for instance the dihedral group $D_8$, for which the angle in the matrix $\rho(a)$ above is $\pi/4$ instead of $2\pi/3$. Hence our choice of the group ${\bf S}_3$ is not arbitrary. In fact, we have long speculated that the surface ${\mathcal S}$ in (\ref{surfaceequ}) contains a higher density of prime triples compared to that of  $\mathbb{Z}^3_+$. Apparently, the symmetry of ${\mathcal S}$ plays a role here, for example if $(z_0,z_1,z_2)$ is a prime triple in ${\mathcal S}$, then so is $(z_2,z_1,z_0)$. Indeed, numerical calculations confirm that among tens of millions of coprime triples in $\mathcal{S}\cap \mathbb{Z}^3_+$, the density of prime triples is greater than that for the set $\mathbb{Z}^3_+$, indicating that $\mathcal{S}$ meets primes more frequently. It is an intriguing question whether this phenomenon holds asymptotically. We will prove that the answer is yes!

First, we give a precise formulation of the problem. Let $\# X$ denote the number of elements in a finite set $X$. For a positive integer $N$, we set $\pi(N) = \#\{~p \leq N~|~p\ \text{is prime}\}$. The Prime Number Theorem \cite[Theorem 6]{HW} asserts that $\pi(N)=\frac{N}{ \log N}(1+o(1))$ when $N$ is large, showing that the density of primes in the set $\{1, \dots, N\}$ is asymptotically equal to $1/\log N$. A simple calculation verifies that if $z$ is a prime triple in ${\mathcal S}$, then $z_0z_1z_2\neq 0$. Hence, we assume no coordinate of $z$ is $0$ and define
\begin{align}\label{DD+}
	D_+(N)=\{z\in {\mathbb{Z}_+^3}|\ 1\leq z_0 \le \frac{6}{5}N, 1\leq z_1 \le N, 1\leq z_2 \le \frac{6}{5}N\}.
\end{align}
The fact that $D_+(N)$ is a cuboid rather than a cube will offer convenience for subsequent study. Define
\begin{align*}  
X_+(N)&= \{z\in D_+(N) \ |~z\ \text{is a prime triple} \}, \\ 
	Y_+(N)&= \{z\in D_+(N) \ |~z\ \text{is a coprime triple}\},
\end{align*}
and set
\[\mathbb{P}_+(N) = \frac{\# X_+(N)}{\# Y_+(N)}.\]
Intuitively speaking, the ratio $\mathbb{P}_+(N)$ is the density of prime triples in the set of coprime triples in ${\mathbb{Z}_+^3}$ bounded by $N$. In Section 2, we will give an estimate of $\mathbb{P}_+(N)$. This will serve as a basis for the subsequent comparison. For the eigensurface $\mathcal{S}$ in \eqref{surfaceequ}, we define 
\[ X_+^\mathcal{S}(N)= {\mathcal S}\cap X_+(N), \hspace{1cm} Y_+^\mathcal{S}(N)= \mathcal{S} \cap  Y_+(N)\]and similarly set
\[\mathbb{P}_+^\mathcal{S}(N) = \frac{\# X_+^\mathcal{S}(N)}{\# Y_+^\mathcal{S}(N)}.\]
In Section 3 and Section 4, we will estimate the sizes of $X_+^\mathcal{S}(N)$
and resp. $Y_+^\mathcal{S}(N)$. It is worth noting that the estimation of $\#Y_+^\mathcal{S}(N)$ is the key to the whole problem. Let $\zeta$ denote the Riemann zeta function. The following is the main theorem of this paper.
\begin{thm}
The following hold for the eigensurface ${\mathcal S}$ of the group ${\bf S}_3$:

\emph{a)} $\liminf \frac{\mathbb{P}_+^\mathcal{S}(N)}{\mathbb{P}_+(N)}\geq \frac{\pi^2}{36(3-\sqrt{6})(\sqrt{2}-1)\zeta(3)}\approx 1.0002$;

\emph{b)} $\limsup \frac{\mathbb{P}_+^\mathcal{S}(N)}{\mathbb{P}_+(N)}\leq \frac{2\sqrt{3}\pi^2}{27\zeta(3)}\approx 1.0534$.
\end{thm}

It is worth mentioning that in the process of proving this theorem, we will establish several facts of independent interest in number theory.  Since this paper concerns the density of prime triples in $\mathbb{Z}^3_+$ and $\mathcal{S}\cap \mathbb{Z}^3_+$, the results presented here can be viewed as higher-dimensional generalizations of the Prime Number Theorem. 


\section{An estimate of $\mathbb{P}_+(N)$}

\par This section gives an estimate of $\mathbb{P}_+(N)$. First, we recall some terminologies and notations which will be used in the sequel. Given two real functions $f$ and $g$ defined on $\mathbb{R}$, with $g$ non-vanishing, we write $f \sim g$ if $\lim_{x\to\infty}f(x)/g(x)= 1$. For $g$ nonnegative, the notation $f = O(g)$ means that there exists a constant $C>0$ such that $\lvert f(x) \rvert \leq Cg(x),\ x\in \mathbb{R}$. Moreover, we write $f = o(g)$ in the case $\lim_{x\to\infty}f(x)/g(x)= 0$. The symbols $ \left  \lfloor\cdot \right  \rfloor$ and $\lceil \cdot \rceil$ denote the floor function and the ceiling function, respectively. The M\"obius function $\mu: {\mathbb Z}_+\to \{0, \pm 1\}$ is defined by
\begin{equation}\label{Mobius}
\mu(d)=\begin{cases}
	1, &\text{if}\ d=1,\\	
	(-1)^k, &\text{if}\ d\ \text{is a product of}\ k\ \text{distinct prime numbers},\\
	0, &\text{if}\  d \ \text{has one or more repeated prime factors}.
\end{cases}
\end{equation}
Furthermore, Nymann \cite{Ny} shows that $k$ positive integers, chosen at random from the set $\{1, \dots, N\}$, are relatively prime with probability $1/\zeta(k)+O(1/N)$ for $k\geq 3$ and probability $1/\zeta(2)+O(\log N/N)$ for $k=2$, where $\zeta$ is the Riemann zeta function. In particular,
\[\zeta(3) = \sum_{n = 1}^{\infty}\frac{1}{n^3} = 1.202056903 \cdots.\]

\begin{thm}\label{propP+}
	The density of prime triples in coprime triples is 
	\[\mathbb{P}_+(N) = \frac{3\zeta(3)}{\log N}(1 + o(1)).\]
\end{thm} 
\pf  We estimate $\# Y_+(N)$ and $\# X_+(N)$ separately. First, it follows from Nymann \cite{Ny} that
\begin{align*}
	\# Y_+(N) &= \sum_{d=1}^{N} \mu(d) \cdot  \left  \lfloor \frac{N}{d} \right  \rfloor \cdot  \left  \lfloor\frac{6N}{5d} \right  \rfloor^2\\
	&= \frac{36}{25\zeta(3)}N^3 + O(N^2)
	= \frac{36}{25\zeta(3)}N^3 + o(N^3).	
\end{align*}

\par To estimate $\# X_+(N)$, one observes that in $D_+(N)$ the total number of triples is 
\[ \left  \lfloor\frac{6}{5}N \right  \rfloor \cdot N \cdot  \left  \lfloor\frac{6}{5}N \right  \rfloor , \] 
and the number of triples with all three elements composite is 
\[ \left(\left  \lfloor\frac{6N}{5} \right  \rfloor- \pi\left(\frac{6N}{5}\right)\right) \left(N - \pi(N) \right) \left(\left  \lfloor\frac{6N}{5} \right  \rfloor- \pi\left(\frac{6N}{5}\right)\right).\] 
Thus, the number of triples with at least one prime is 
\begin{align*}
	&\mathrm{\hphantom{={}}}  \left  \lfloor\frac{6}{5}N \right  \rfloor \cdot N \cdot  \left  \lfloor\frac{6}{5}N \right  \rfloor - \left(\left  \lfloor\frac{6N}{5} \right  \rfloor- \pi\left(\frac{6N}{5}\right)\right) \left(N - \pi(N) \right) \left(\left  \lfloor\frac{6N}{5} \right  \rfloor- \pi\left(\frac{6N}{5}\right)\right) \\
	&= \frac{108N^3}{25\log N} + o(\frac{N^3}{\log N}).
\end{align*}
We then subtract the number of such triples which are not relatively prime. This occurs only if one of $z_0, z_1, z_2$ is a prime $p$ and the other two are integer multiples of $p$. The number of such triples is 
\begin{align*}
&\sum_{p \leq N}\left( \left  \lfloor\frac{ \left  \lfloor\frac{6}{5}N \right  \rfloor}{p} \right  \rfloor \right)^2 + 2\sum_{p \leq \frac{6}{5}N}\left( \left  \lfloor\frac{ \left  \lfloor\frac{6}{5}N \right  \rfloor}{p} \right  \rfloor  \left  \lfloor\frac{N}{p} \right  \rfloor \right) \\
\leq{} &3\sum_{p \leq \frac{6}{5}N }\left(\frac{6N}{5p}\right)^2 \leq O(N^2) = o\left(\frac{N^3}{\log N}\right).
\end{align*}
Therefore, 
\[\# X_+(N) = \frac{108N^3}{25\log N} + o\left(\frac{N^3}{\log N}\right), \]	
and it follows that
\[\mathbb{P}_+(N) = \frac{\# X_+(N)}{\# Y_+(N)} = \frac{3\zeta(3)}{\log N}(1 +o(1)).\]  \zb\\

\section{An estimate of $\# X_+^\mathcal{S}(N)$}

This section gives an estimate of $\#X_+^\mathcal{S}(N)$, where we recall that
\begin{align*}  
	X_+^\mathcal{S}(N)= \{z\in \mathcal{S} \cap D_+(N) \ |~z\ \text{is a prime triple}\}.
\end{align*}

\begin{thm}\label{X+} 
	$ \# X_+^\mathcal{S}(N)= \frac{2N}{\log N} + o\left(\frac{N}{\log N}\right).$
\end{thm}

To prove Theorem \ref{X+}, we need to parameterize the eigensurface 
$\mathcal{S}$. Recall that $\mathcal{S}$ is determined by the Diophantine equation 
\begin{align}\label{0231}
	z_0^{2} - z_1^{2} + z_2^{2} - z_0z_2= 0. 
\end{align}
Rewrite equation \eqref{0231} as
\begin{align}\label{0232}
	z_1^{2}= z_0^{2} + z_2^{2} - z_0z_2= z_0^{2} + z_2^{2} - 2z_0z_2\cos \frac{\pi}{3}.
\end{align}
By the Law of Cosines, a positive integer solution of equation \eqref{0232} corresponds to an integer-sided triangle with a $\frac{\pi}{3}$ angle, whose opposite side is $z_1$. If $z_0=z_2$, then $z_0=z_1=z_2$. If $z_0\neq z_2$, Gilder \cite{Gi} proved the following fact about the non-equilateral integer-sided triangles with a $\frac{\pi}{3}$ angle.

\begin{lem}\label{IT_60}	
	A non-equilateral integer-sided triangle has an angle of $\frac{\pi}{3}$ if and only if it is similar to precisely one member of the class of triangles whose sides are $m^2 + mn + n^2, m^2 + 2mn$ and either $n^2 + 2mn$ or $m^2 - n^2$, where $m$ and $n$ are positive coprime integers with $m > n$ and $m \not\equiv n \pmod 3$. 
	\par In addition, for each triangle in this class, the greatest common divisor of the three sides is $1$.
\end{lem}

Lemma \ref{IT_60} offers a parameterization of coprime solutions of equation \eqref{0231}. Set
\[\Omega=\{(m, n)\in \mathbb{Z}_+^2|\ \gcd(m,n)=1,\  m > n \ \text{and}\  m \not\equiv n \pmod 3\}\]
and define the following four maps
\begin{align*}
	&\phi_1: \mathbb{Z}_+^2 \to \mathbb{Z}_+^3, (m, n) \mapsto (m^2 + 2mn, m^2 + mn + n^2, n^2 + 2mn), \\
	&\phi_2: \mathbb{Z}_+^2 \to \mathbb{Z}_+^3, (m, n) \mapsto (n^2 +2mn, m^2 + mn + n^2, m^2 + 2mn), \\
	&\phi_3: \mathbb{Z}_+^2 \to \mathbb{Z}_+^3, (m, n) \mapsto (m^2 + 2mn, m^2 + mn + n^2, m^2 - n^2), \\
	&\phi_4: \mathbb{Z}_+^2 \to \mathbb{Z}_+^3, (m, n) \mapsto (m^2 - n^2, m^2 + mn + n^2, m^2 + 2mn).
\end{align*}
It is easy to verify that all four maps are injective. If $m> n$, the ranges of those maps are pairwise disjoint. Hence we have the following lemma. 	

\begin{lem}\label{PIS}
	Given any $z\in \mathbb{Z}_+^3$ with $\gcd(z)=1$, it is a solution of equation \eqref{0231} if and only if it belongs to one of the following five disjoint sets:
	\begin{flalign*}
		&\emph{(i)}\ \Delta_0=\{(1, 1, 1)\},&\\
		&\emph{(ii)}\ \Delta_1 = \{\phi_1(m, n)\ |\ (m, n) \in \Omega\},&\\
		&\emph{(iii)}\ \Delta_2 = \{\phi_2(m, n)\ |\ (m, n) \in \Omega\},&\\
		&\emph{(iv)}\ \Delta_3 = \{\phi_3(m, n)\ |\ (m, n) \in \Omega\},&\\
		&\emph{(v)}\ \Delta_4= \{\phi_4(m, n)\ |\ (m, n) \in \Omega\}.&
	\end{flalign*}		
\end{lem}

Thus, counting $X_+^\mathcal{S}(N)$ amounts to computing the number of prime triples in each $\Delta_i \cap D_+(N)$. Observe that the symmetry of ${\mathcal S}$ in $z_0$ and $z_2$ implies 
\begin{equation}\label{symm}
\begin{aligned}
 &\# (\Delta_1 \cap D_+(N))=\# (\Delta_2 \cap D_+(N)), \\ 
 &\# (\Delta_3 \cap D_+(N))=\# (\Delta_4 \cap D_+(N)).
\end{aligned}
 \end{equation}
The next lemma justifies the use of cuboid $D_+(N)$ (see (\ref{DD+})) instead of a cube. 
\begin{lem}\label{D+}
For any $z=(z_0, z_1, z_2)\in \Delta_i$, $1\leq i\leq 4$, it holds that $z\in D_+(N)$ if and only if $z_1 \leq N$.	
\end{lem}	

\pf It suffices to prove the sufficiency. We will show that for any $z\in \Delta_1$, when $z_1 \leq N$, we have $z\in D_+(N)$. The proofs for $\Delta_i$, $i=2, 3, 4$ are similar.
\par For any $(z_0, z_1, z_2)\in \Delta_1$, we have $z_0=m^2 + 2mn$, $z_1=m^2 + mn + n^2$, $z_2=n^2 + 2mn$ for some $(m,n)\in\Omega$. Note that $m>n$, hence if $z_1 \leq N$, we easily get  
$z_2<z_1\leq N < \frac{6}{5}N$. Since 
\[\frac{z_0}{z_1}=\frac{m^2 +2mn}{m^2 +mn +n^2}=\frac{\left(\frac{m}{n}\right)^2 +2\left(\frac{m}{n}\right)}{(\frac{m}{n})^2 +\frac{m}{n}+1}\leq \frac{2}{\sqrt{3}}< \frac{6}{5}, \]
we have that $z_0<\frac{6}{5}z_1 \leq \frac{6}{5}N$, showing $z\in D_+(N)$. \zb\\

Two other functions are needed to prove Theorem \ref{X+}. The first is the function
\[\pi(N; q, a) := \#\{p \leq  N \ | \ p\ \text{is prime}, p\equiv a \pmod q\},\]
which counts the number of primes $p \equiv a \pmod q$ up to $N$; and the second is Euler's totient function $\varphi(q)$, which counts how many numbers from $1$ to $q$ are relatively prime to $q$. Landau's theorem gives an estimate of $\pi(N; q, a)$ (see \cite[Section 120]{La} and \cite[Theorem 7.25]{LeV}). 

\begin{lem}\label{Lan}
	Let $a, q > 0$ be relatively prime integers. Then 
	\[\pi(N; q, a) = \frac{1}{\varphi(q)} \frac{N}{\log N} \big(1 + o(1)\big).\]
\end{lem}

For instance, since $\varphi(6)=2$, asymptotically we have \[\pi(N; 6, 1)=\pi(N; 6, 5)=\frac{N}{2\log N}(1+o(1)),\] which is a half of the number of primes up to $N$. Now we are ready to prove Theorem \ref{X+}.

\noindent\textbf{Proof of Theorem \ref{X+}.} 
Recall that the set $X_+^\mathcal{S}(N)$ consists of prime triples $z\in \mathcal{S}\cap D_+(N)$. Obviously, the set $\Delta_0$ does not contain such triples. Thus, in view of Lemma \ref{PIS} and (\ref{symm}), it suffices to count such triples in each $\Delta_i\cap D_+(N)$, $i=1, 3$.

 For $z \in \Delta_1\cap D_+(N)$, since $z_0 = m^2 + 2mn$ can not be prime, we consider the following two subsets of $\Delta_1$. Let 
\begin{align*}
	\Delta_1^1=\{z\in\Delta_1 |\ z_1 \text{ is prime}\},\hspace{1cm}\Delta_1^2=\{z\in\Delta_1|\ z_2 \text{ is prime}\}.
\end{align*}
We will count $\Delta_1^1\cap D_+(N)$ and $\Delta_1^2\cap D_+(N)$ separately. 

For $z\in \Delta_1^1$, Lemma \ref{D+} indicates that $z\in \Delta_1^1 \cap D_+(N)$ if and only if $z_1\leq N$.
Since $z_1=m^2+mn+n^2$, we need to count the number of coprime pairs $(m,n)\in \Omega$ such that $m^2+mn+n^2$ is a prime not exceeding $N$. Nair \cite{Na} proved that for each prime $p=m^2 + mn + n^2$, the corresponding coprime pair $(m,n)$ with $m>n$ is unique, and hence corresponds to a point in $\Delta_1^1$. Further, Gilder \cite{Gi} observed that a prime $p\neq 3$ is of the form $m^2 + mn + n^2$ for some coprime pair $(m,n)$ if and only if $p=6k + 1$, for some $k \in {\mathbb{Z_+}}$. Thus, Lemma \ref{Lan} implies
\[\# (\Delta_1^1\cap D_+(N))= \frac{1}{2}\frac{N}{\log N}(1 + o(1)).\]

For $z\in\Delta_1^2$, $z_2 = n^2 + 2mn$ can be prime only if $n = 1$. Hence
\[\#(\Delta_1^2\cap D_+(N)) \leq\#\{\phi_1(m, 1)\in D_+(N)\ |\ (m, 1) \in \Omega\}.\]
Since $\phi_1(m, 1)=(m^2 + 2m, m^2 + m + 1, 2m + 1)$, by Lemma \ref{D+}, $\phi_1(m, 1)\in D_+(N)$ if and only if $m^2+m+1\le N$, i.e., $m \leq \frac{\sqrt{4N - 3} - 1}{2}$. Therefore,
\[\#(\Delta_1^2\cap D_+(N)) \leq\left \lfloor\sqrt{4N-3} \right  \rfloor = o\left(\frac{N}{\log N}\right),\]
and consequently,
\[\#(\Delta_1^1\cap \Delta_1^2\cap D_+(N) )\leq\#(\Delta_1^2\cap D_+(N)) \leq o\left(\frac{N}{\log N}\right).\]
Hence, we have
\begin{align*}
	\#(\Delta_1 \cap D_+(N))
	&=\#(\Delta_1^1\cap D_+(N))+\#(\Delta_1^2\cap D_+(N))-\#(\Delta_1^1\cap \Delta_1^2\cap D_+(N) )\\ &=\frac{1}{2}\frac{N}{\log N}(1 + o(1)).
\end{align*}

The estimation for $\Delta_3\cap D_+(N)$ is similar. In summary, we have
\[ \# X_+^\mathcal{S}(N)= \frac{2N}{\log N} + o\left(\frac{N}{\log N}\right)\]    
as claimed by Theorem \ref{X+}.\zb\\

\section{An estimate of $\# Y_+^\mathcal{S}(N)$}

\par This section gives an estimate of $\# Y_+^\mathcal{S}(N)$. This is the critical step for establishing the main theorem. First, recall that 
\begin{align*}  
	Y_+^\mathcal{S}(N)&= \{z\in \mathcal{S} \cap D_+(N) \ |\gcd z= 1\}.
\end{align*}

\begin{thm}\label{Y+}
	The following inequalities hold for $Y_+^\mathcal{S}(N)$:
	\[\frac{3\sqrt{3}}{\pi^2}N + o(N)< \# Y_+^\mathcal{S}(N) <\frac{24(3-\sqrt{6})(\sqrt{2}-1)}{\pi^2}N + o(N).\] 
\end{thm}

Like the case for $X_+^\mathcal{S}(N)$, counting $Y_+^\mathcal{S}(N)$ amounts to counting $\Delta_i\cap D_+(N)$ for each $i=1, 3$, and it is sufficient to do so for $i=1$.

\par
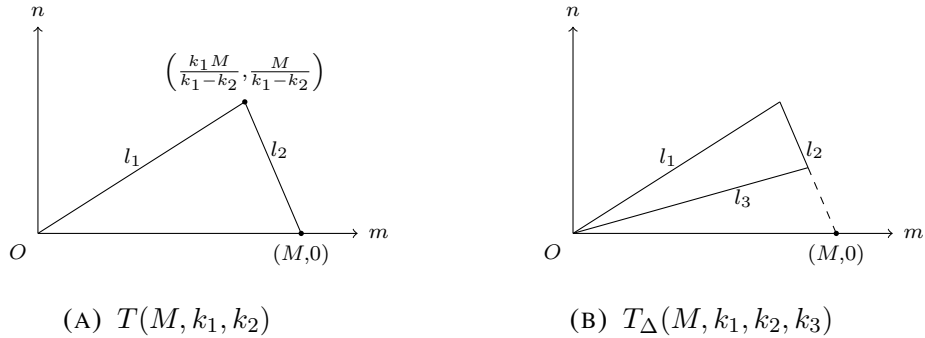
\begin{figure}[ht]
	\centering
	\begin{subfigure}[b]{0.32\textwidth}
		\centering
		\begin{tikzpicture}[scale=2.5]
			\draw[->, thin] (0,0) -- (1.7,0) node[right] {$\scriptstyle m$};
			\draw[->, thin] (0,0) -- (0,1.1) node[above] {$\scriptstyle n$};
			\draw[thin] (0,0) node[below left] {$\scriptstyle O$};
			
			\draw[black, thin] (0,0) -- (1.1,0.7);
			\node[black, above right] at (0.4,0.3) {$\scriptstyle l_1$};
			
			\draw[black, thin] (1.1,0.7) -- (1.4,0);
			\node[black, above right] at (1.18,0.34) {$\scriptstyle l_2$};
			
			\fill[black] (1.1,0.7) circle (0.4pt) node[above] {$\scriptstyle \left(\frac{k_1M}{k_1 - k_2}, \frac{M}{k_1 - k_2}\right)$};
			\fill[black] (1.4,0) circle (0.4pt) node[below] {$\scriptstyle (M, 0)$};
		\end{tikzpicture}
		\caption{~$T(M, k_1, k_2)$}
		\label{subfig:11}
	\end{subfigure}
	\qquad\qquad\qquad
	\begin{subfigure}[b]{0.32\textwidth}
		\centering
		\begin{tikzpicture}[scale=2.5]
			\draw[->, thin] (0,0) -- (1.7,0) node[right] {$\scriptstyle m$};
			\draw[->, thin] (0,0) -- (0, 1.1) node[above] {$\scriptstyle n$};
			\draw[thick] (0,0) node[below left] {$\scriptstyle O$};
			
			\draw[black, thin] (0,0) -- (1.1,0.7);
			\node[black, above right] at (0.4,0.3) {$\scriptstyle l_1$};
			
			\draw[black, thin] (1.1,0.7) -- (1.25,0.35);
			\node[black, above right] at (1.18,0.34) {$\scriptstyle l_2$};
			
			\draw[black, dashed] (1.25,0.35) -- (1.4,0);
			
			\draw[black, thin] (0,0) -- (1.25,0.35);
			\node[black, above right] at (0.8,0.08) {$\scriptstyle l_3$};
			

			
			
			\fill[black] (1.4,0) circle (0.4pt) node[below] {$\scriptstyle (M, 0)$};
		\end{tikzpicture}
		\caption{~$T_\Delta(M, k_1, k_2, k_3)$}
		\label{subfig:12}
	\end{subfigure}
	\caption{Diagram for the triangular regions (with $k_2<0$)}\label{fig:1}
\end{figure}

 For $M>0$, $k_1>0$ and $k_1>k_2$, let $T(M, k_1, k_2)$ be the  triangular region within the first quadrant of the plane, enclosed by the $m-$axis, the lines $l_1: m = k_1n$ and $l_2: m = k_2n + M$, as shown in \textsc{Figure} \ref{fig:1} (A). Define 
\begin{align*}
	&C(M, k_1, k_2) = \#\{(m,n)\in T(M, k_1, k_2)~|~\gcd(m,n)=1\},\\
	&\hat{C}(M, k_1, k_2) = \#\{(m,n)\in T(M, k_1, k_2)~|~\gcd(m,n)=1,~m \not\equiv n \pmod 3 \}.
\end{align*}
We will estimate $C(M, k_1, k_2)$ and $\hat{C}(M, k_1, k_2)$ below. The estimation of the latter is subtle, but it is a key step to the computation of $\# Y_+^\mathcal{S}(N)$.

\begin{lem}\label{TRI} 
The following hold asymptotically for large $M$:	
	\begin{flalign*}
	&\emph{(i)}\ C(M, k_1, k_2) = \frac{3}{\pi^2}\frac{M^2}{k_1- k_2} + o(M^2),& \\
	&\emph{(ii)}\ \hat{C}(M, k_1, k_2)=\frac{3}{4}C(M, k_1, k_2) + o(M^2).
\end{flalign*}
\end{lem}
\par In order to prove Lemma \ref*{TRI}, we need some basic facts about the sums of the Euler's totient function $\varphi(n)$ and the function $d(n)$ which counts the number of positive divisors of $n$. They are summarized as follows (see \cite[Section 3]{Ap} and \cite[Theorem 320]{HW}).

\begin{lem}\label{EUL} For sufficiently large $N$, the following
	 hold asymptotically:
	\begin{flalign*}
		&\emph{(i)}\ \sum_{n=1}^{N}\varphi(n) = \frac{3}{\pi^2}N^2 + O(N \log N),& \\
		&\emph{(ii)}\ \sum_{n=1}^{N}\frac{\varphi(n)}{n} = \frac{6}{\pi^2}N + O(\log N),\\
		&\emph{(iii)}\ \sum_{n=1}^{N}d(n) = N \log N +O(N).
	\end{flalign*}
\end{lem}

\par We are now in position to prove Lemma \ref{TRI}.

\noindent\textbf{Proof of Lemma \ref{TRI}.} Given a positive integer $n$, let $f(N, n)$ denote the number of integers between $1$ and $N$ that are relatively prime to $n$. Obviously, $f(n, n)$ is Euler's
totient function $\varphi(n)$. T\'oth \cite{To} shows that
\[f(N,n) = \sum_{d \mid n}\mu(d)\left  \lfloor\frac{N}{d} \right  \rfloor,\]where $\mu(d)$ is the M\"obius function (\ref{Mobius}). Let $\omega(n)$ be the number of distinct prime factors of $n$. Then, as a consequence of \cite[Section 16.3]{HW}, we have
\[f(N,n) = \sum_{d \mid n}\mu(d)\left  \lfloor\frac{N}{d} \right  \rfloor= \frac{\varphi(n)}{n}N + O\left(2^{\omega(n)}\right).\]

\par  Thus in $T(M, k_1, k_2)$, for a fixed $n$, the number of $m \in [k_1n, k_2n + M] $ that are relatively prime to $n$ is 
\[f(\left  \lfloor k_2n + M  \right  \rfloor, n)-f(\lceil k_1n \rceil,n) = \frac{\varphi(n)}{n}( \left  \lfloor k_2n + M  \right  \rfloor - \lceil k_1n \rceil) + O\left(2^{\omega(n)}\right).\]
Set $h =  \left  \lfloor\frac{M}{k_1 - k_2} \right  \rfloor$. Then the number of coprime pairs $(m, n)$ in $T(M, k_1, k_2)$ is 
\begin{align*}
	C(M, k_1, k_2) &= \sum_{n=1}^{h}\frac{\varphi(n)}{n}( \left  \lfloor k_2n + M  \right  \rfloor - \lceil k_1n \rceil) + O\bigg(\sum_{n=1}^{h}2^{\omega(n)}\bigg)\\
	&= \sum_{n=1}^{h}\frac{\varphi(n)}{n}(k_2n + M - k_1n + O(1)) + O\bigg(\sum_{n=1}^{h}2^{\omega(n)}\bigg)\\
	&= (k_2- k_1)\sum_{n=1}^{h}\varphi(n) + M\sum_{n=1}^{h}\frac{\varphi(n)}{n} + O\bigg(\sum_{n=1}^{h}\frac{\varphi(n)}{n}\bigg) + O\bigg(\sum_{n=1}^{h}2^{\omega(n)}\bigg).
\end{align*}
From Lemma \ref{EUL} (i) and (ii), we have 
\begin{equation}\label{0333}
	\begin{aligned}
		C(M, k_1, k_2) 
		={}& \frac{3}{\pi^{2}}(k_2- k_1)h^2 + \frac{6}{\pi^{2}}Mh \\& + O(h\log h)+ O(M\log h) + O(h) + O\bigg(\sum_{n=1}^{h}2^{\omega(n)}\bigg).
	\end{aligned}
\end{equation}
Substituting $h = \left \lfloor\frac{M}{k_1 - k_2} \right  \rfloor = \frac{M}{k_1-k_2} + O(1)$ into \eqref{0333}, we have 
\[
C(M, k_1, k_2)
=\frac{3}{\pi^2}\frac{M^2}{k_1 - k_2} +o(M^2) + O\bigg(\sum_{n=1}^{h}2^{\omega(n)}\bigg). 
\]
\par Now we estimate $\displaystyle O\bigg(\sum_{n=1}^{h}2^{\omega(n)}\bigg)$.
For a given $n = p_1^{t_1}p_2^{t_2} \dots p_k^{t_k}$, where $p_1, p_2, \dots, p_k$ are different primes and $t_1, t_2, \dots, t_k \in \mathbb{Z_+}$, we have that \[\omega(n) = k, \ d(n) = (t_1 + 1)(t_2 + 1) \cdots (t_k + 1).\]
Since $2 \leq t_i + 1$ for each $i$, we have $2^{\omega(n)} \leq d(n)$.
Hence by Lemma \ref{EUL} (iii),
\[\sum_{n=1}^{h}2^{\omega(n)} \leq \sum_{n=1}^{h}d(n) = h\log h +O(h),\]
which means $\displaystyle O\bigg(\sum_{n=1}^{h}2^{\omega(n)}\bigg) = o(M^2)$.
Therefore, 
\[C(M, k_1, k_2) = \frac{3}{\pi^2}\frac{M^2}{k_1- k_2} + o(M^2).\] 

\par (ii) Next, we set \begin{align*}	
	C'(M, k_1, k_2) = \#\{(m,n)\in T(M, k_1, k_2) | \gcd(m,n)=1, m \equiv n \pmod 3 \}
\end{align*}
and show that  
\[C'(M, k_1, k_2) =  \frac{1}{4}C(M, k_1, k_2) + o(M^2).\]
The claim in Lemma \ref{TRI} (ii) will then immediately follows. 

We need to count the number of coprime pairs $(m,n)$ in $T(M, k_1, k_2)$ satisfying $m \equiv n \pmod 3$. Let $n$ be fixed, then $m \in [k_1n, k_2n + M]$. Since $m \equiv n \pmod 3$, we can write $m = n + 3t$, where $t$ is an integer in $\left[\frac{(k_1 -1)n}{3}, \frac{(k_2-1)n + M}{3}\right]$. Due to the fact $\gcd(m, n)= \gcd(n+3t, n) = \gcd(3t, n)$, the number of coprime pairs $(m, n)$ in $T(M, k_1, k_2)$ is equal to that of $t$ in $\left[\frac{(k_1 -1)n}{3}, \frac{(k_2-1)n + M}{3}\right]$ for which $\gcd(3t, n)=1$. Note that $\gcd(3t, n)=1$ if and only if $3 \nmid n$ and $\gcd(t, n)=1$. Hence $C'(M, k_1, k_2) $ is equal to the number of coprime pairs $(t,n)$ satisfying $n\in [1,h]$, $3\nmid n$ and $t\in\left[\frac{(k_1 -1)n}{3}, \frac{(k_2-1)n + M}{3}\right]$. Thus,
\begin{align*}
	C'(M, k_1, k_2)= \sum_{\substack{n=1 \\ 3\nmid n}}^{h}\frac{\varphi(n)}{n}\left( \left  \lfloor\frac{(k_2-1)n + M}{3} \right  \rfloor- \left \lceil \frac{(k_1 -1)n}{3} \right \rceil \right) + O\bigg(\sum_{\substack{n=1 \\ 3\nmid n}}^{h}2^{\omega(n)}\bigg).	
\end{align*} 
By the same reasoning as that in part (i), we have 
\begin{align*}
	C'(M, k_1, k_2)
	= \left(\frac{k_2 - k_1}{3}\right)\sum_{\substack{n=1 \\ 3\nmid n}}^{h}\varphi(n) + \frac{M}{3}\sum_{\substack{n=1 \\ 3\nmid n}}^{h}\frac{\varphi(n)}{n} + o(M^2).
\end{align*}
The sums
$\displaystyle\sum_{\substack{n=1 \\ 3\nmid n}}^{h}\varphi(n)\  \text{and}\  \sum_{\substack{n=1 \\ 3\nmid n}}^{h}\frac{\varphi(n)}{n}$ need to be estimated separately.

\par Obviously,
\begin{align*}
	\sum_{\substack{n=1 \\ 3\nmid n}}^{h}\varphi(n) = \sum_{n=1}^{h}\varphi(n) - \sum_{\substack{n=1 \\ 3\mid n}}^{h}\varphi(n).
\end{align*}
For $N\in \mathbb{Z}_+$, set $ \Phi(N)=\sum_{n=1}^{N}\varphi(n)$ and $\Phi(N,3)=\sum_{\substack{n=1 \\ 3\mid n}}^{N}\varphi(n)$. The estimate of $\Phi(N)$ follows from Lemma \ref{EUL} (i), so it remains to estimate $\Phi(N, 3)$. The estimation procedure follows  Lehmer's computation in \cite[Lemma 2]{Le}.

First, for $j \in \mathbb{Z_+}$ we will prove
\begin{equation}\label{varphi0}
	\begin{aligned}\varphi(3j) = \begin{cases} 
			2\varphi(j), &\text{if}\ 3 \nmid j,\\
			3\varphi(j), &\text{if}\ 3 \mid j.
		\end{cases}
	\end{aligned}
\end{equation}
It is not hard to see that if $\gcd(a, b)=1$, then  $\varphi(ab)=\varphi(a)\varphi(b)$ \cite[Theorem 60]{HW}. Hence, if $3 \nmid j$, then $\varphi(3j) = \varphi(3)\varphi (j)= 2\varphi(j)$.
If $3 \mid j$, there exist $q, t \in \mathbb{Z_+}$ such that $j = 3^qt$ and $\gcd(3, t)=1$. Observe that $\varphi(3^{q+1}) = 3^{q+1} - 3^q$, and consequently,
\[\varphi(3j) = \varphi(3^{q+1}t)= \varphi(3^{q+1})\varphi(t) = (3^{q+1} - 3^q)\varphi(t).\] 
Similarly, we have $\varphi(j)= \varphi(3^qt)= (3^{q} - 3^{q-1})\varphi(t)$ and hence $\varphi(3j)=3\varphi(j)$.
This proves equation \eqref{varphi0}.
Next, we will establish the following iterative formula. For $N\in \mathbb{Z}_+$, it holds that
\begin{align} \label{varphi1}
	\Phi(N,3) = 2\Phi(\left  \lfloor  N/3 \right \rfloor) + \Phi(\left  \lfloor  N/3 \right \rfloor, 3).
\end{align}
In fact, since
\begin{align*} 
	\Phi(N,3) = \sum_{\substack{n=1 \\ 3\mid n}}^{N}\varphi(n) = \sum_{j=1}^{\left  \lfloor  N/3 \right \rfloor}\varphi(3j)=\sum_{\substack{j = 1 \\ 3 \nmid j}}^{ \left  \lfloor N/3 \right \rfloor}\varphi(3j) + \sum_{\substack{j = 1 \\ 3 \mid j}}^{ \left  \lfloor N/3 \right \rfloor}\varphi(3j),
\end{align*}
by equation \eqref{varphi0}, we have
\begin{align*} 
	\Phi(N,3) &= 2\sum_{\substack{j = 1 \\ 3 \nmid j}}^{ \left  \lfloor N/3 \right \rfloor}\varphi(j) + 3\sum_{\substack{j = 1 \\ 3 \mid j}}^{ \left  \lfloor N/3 \right \rfloor}\varphi(j)\\
	&= 2\sum_{j= 1}^{ \left \lfloor N/3 \right \rfloor}\varphi(j) - 2\sum_{\substack{j = 1 \\ 3 \mid j}}^{ \left  \lfloor N/3 \right \rfloor}\varphi(j) + 3\sum_{\substack{j = 1 \\ 3 \mid j}}^{ \left  \lfloor N/3 \right \rfloor}\varphi(j)\\
	&= 2\Phi(\left  \lfloor  N/3 \right \rfloor) + \Phi(\left  \lfloor  N/3 \right \rfloor, 3).
\end{align*}
Now we apply equation \eqref{varphi1} repeatedly to estimate $\Phi(h,3)$. Note that for $a$, $b \in \mathbb{Z_+}$,  we have $\left \lfloor \left  \lfloor h/a \right \rfloor /b  \right \rfloor = \lfloor h/ab \rfloor$ \cite[Section 3.2]{GKP}. It follows that 
\begin{align*} 
	\Phi(h,3) &= 2\Phi(\left  \lfloor  h/3 \right \rfloor) + \Phi(\left  \lfloor  h/3 \right \rfloor, 3)\\
	&= 2\Phi(\left  \lfloor  h/3 \right \rfloor) + \left( 2\Phi(\left  \lfloor  h/3^2 \right \rfloor) + \Phi(\left  \lfloor  h/3^2 \right \rfloor, 3)  \right)\\
	&=2\Phi(\lfloor  h/3 \rfloor) + 2\Phi(\lfloor  h/3^2 \rfloor) + \cdots + 2\Phi(\lfloor h/3^{\lfloor \log_3 h \rfloor} \rfloor).
\end{align*}
Noting that the above iterative process has only $\lfloor \log_3 h \rfloor$ steps. Substituting Lemma \ref{EUL} (i) into the above formula, we obtain
\begin{align*} 
	\Phi(h,3) =\frac{3}{4\pi^2}h^2 + o(h^2)
\end{align*}
and consequently
\begin{align*}
	\sum_{\substack{n=1 \\ 3\nmid n}}^{h}\varphi(n) = \Phi(h)-\Phi(h,3)= \frac{9}{4\pi^2}h^2 +o(h^2).
\end{align*}

\par Next, we estimate the sum $ \sum_{\substack{n=1 \\ 3\nmid n}}^{h}\frac{\varphi(n)}{n}$. Likewise, we write it as $\sum_{n=1}^{h}\frac{\varphi(n)}{n} - \sum_{\substack{n=1 \\ 3\mid n}}^{h}\frac{\varphi(n)}{n}$, 
and for $N\in \mathbb{Z}_+$, set $ \Phi_1(N)=\sum_{n=1}^{N}\frac{\varphi(n)}{n}$ and $ \Phi_1(N,3)=\sum_{\substack{n=1 \\ 3\mid n}}^{N}\frac{\varphi(n)}{n}$. The key is to compute $\Phi_1(N,3)$. Similar to the preceding procedure, we shall establish the following iterative formula
\begin{align} \label{varphi2}
	\Phi_1(N,3) = \frac{2}{3}\Phi_1(\left  \lfloor  N/3 \right \rfloor) + \frac{1}{3}\Phi_1(\left  \lfloor  N/3 \right \rfloor, 3).
\end{align}
We apply equation \eqref{varphi2} repeatedly to estimate $\Phi_1(h,3)$. It follows that 
\begin{align*} 
	\Phi_1(h,3) &= \frac{2}{3}\Phi_1(\left  \lfloor  h/3 \right \rfloor) + \frac{1}{3}\Phi_1(\left  \lfloor  h/3 \right \rfloor, 3)\\
	&= \frac{2}{3}\Phi_1(\left  \lfloor  h/3 \right \rfloor) + \frac{1}{3}\left( \frac{2}{3}\Phi_1(\left  \lfloor  h/3^2 \right \rfloor) + \frac{1}{3}\Phi_1(\left  \lfloor  h/3^2 \right \rfloor, 3)  \right)\\
	&=\frac{2}{3}\Phi_1(\lfloor  h/3 \rfloor) + \frac{1}{3}\cdot\frac{2}{3}\Phi_1(\lfloor  h/3^2 \rfloor) + \cdots \\ &+\frac{1}{3^{\lfloor \log_3 h \rfloor-1}}\cdot\frac{2}{3}\Phi_1(\lfloor h/3^{\lfloor \log_3 h \rfloor} \rfloor).
\end{align*}
Substituting Lemma \ref{EUL} (ii) into the above formula, we have 
\begin{align*} 
	\Phi_1(h,3) =\frac{3}{2\pi^2}h + o(h),
\end{align*}
and hence
\begin{align*}
	\sum_{\substack{n=1 \\ 3\nmid n}}^{h}\frac{\varphi(n)}{n} = \displaystyle \Phi_1(h)-\Phi_1(h,3)= \frac{9}{2\pi^2}h +o(h).
\end{align*}
Therefore,
\begin{align*}
	&C'(M, k_1, k_2)
	=\left(\frac{k_2 - k_1}{3}\right)\sum_{\substack{n=1 \\ 3\nmid n}}^{h}\varphi(n) + \frac{M}{3}\sum_{\substack{n=1 \\ 3\nmid n}}^{h}\frac{\varphi(n)}{n} + o(M^2)\\
	={}& \left(\frac{k_2 - k_1}{3}\right)\frac{9}{4\pi^2}h^2 + \frac{M}{3}\cdot\frac{9}{2\pi^2}h + \left(\frac{k_2 - k_1}{3}\right)o(h^2) +\frac{M}{3}o(h)+ o(M^2).
\end{align*}
Substituting $h = \left \lfloor\frac{M}{k_1 - k_2} \right  \rfloor = \frac{M}{k_1-k_2} + O(1)$ into the above formula, we have
\[C'(M, k_1, k_2) =\frac{3}{4\pi^2}\frac{M^2}{k_1 - k_2} + o(M^2),\]
and consequently, 
\[C'(M, k_1, k_2)=\frac{1}{4}C(M, k_1, k_2) +o(M^2),\] 
giving the desired result. \zb\\

Lemma \ref{TRI} can be easily generalized to the following situation. For $M>0$, $k_3>k_1>0$ and $k_1>k_2$, let $T_\Delta(M, k_1, k_2, k_3)$ be the triangular region within the first quadrant of the plane, enclosed by $l_1: m = k_1n$, $l_2: m = k_2n + M$ and $l_3: m=k_3n$, as shown in \textsc{Figure} \ref{fig:1} (B). Let $C_\Delta(M, k_1, k_2, k_3)$ be the number of coprime pairs $(m, n)$ in $T_\Delta(M, k_1, k_2, k_3)$, and $\hat{C}_\Delta(M, k_1, k_2, k_3)$ be the number of coprime pairs $(m, n)$ in $T_\Delta(M, k_1, k_2, k_3)$ such that $m \not\equiv n \pmod 3$. Then we have the following corollary.

\begin{corr}\label{TDelta} 
The following estimates hold:
	\begin{flalign*}
		&\emph{(i)}\ C_\Delta(M, k_1, k_2, k_3) =\frac{3}{\pi^2}\left(\frac{1}{k_1- k_2} - \frac{1}{k_3- k_2}\right)M^2 + o(M^2),& \\
		&\emph{(ii)}\ \hat{C}_\Delta(M, k_1, k_2, k_3)=\frac{3}{4}C_\Delta(M, k_1, k_2, k_3) + o(M^2).
	\end{flalign*}
\end{corr}

\par We remark here that the area of the triangle $T(M, k_1, k_2)$ is 
\[S(M, k_1, k_2)=\frac{1}{2}\frac{M^2}{k_1-k_2}.\] Thus, Lemma \ref{TRI} (i) can be rewritten as
\[C(M, k_1, k_2) = \frac{6}{\pi^2}S(M, k_1, k_2) + o(M^2).\] 
Corollary \ref{TDelta} (i) can also be rewritten as
\begin{align}\label{area}
	C_\Delta(M, k_1, k_2, k_3) =\frac{6}{\pi^2}S_\Delta(M, k_1, k_2, k_3) + o(M^2),
\end{align}	
where $S_\Delta(M, k_1, k_2, k_3)$ is the area of the triangle $T_\Delta(M, k_1, k_2, k_3)$. Formula \eqref{area} is what we actually use in the following computation. 

It is a classical result that, the density of coprime pairs in a ``reasonable" planar region is approximately $1/\zeta(2)=6/\pi^2$ (see \cite{Ro} and \cite[Theorem 459]{HW}). Lemma \ref{TRI} (i) offers an elementary proof of this fact for any triangular region. The key point is Lemma \ref{TRI} (ii), which shows that the density of coprime pairs $(m,n)$ such that $m \not\equiv n \pmod 3$ in the triangular region is $\frac{3}{4}\frac{6}{\pi^2}$. In fact, using the same method as that in the proof of Lemma \ref{TRI} (ii), we arrive at the following variation of Nymann's theorem for $k=2$. 


\begin{corr}
Let $p$ be a prime. The probability that a pair of positive integers $(m, n)$ chosen at random is coprime and fulfilling $m \not\equiv n \pmod p$ is $\frac{p}{p+1}\frac{6}{\pi^2}$.
\end{corr}

We conjecture that this fact remains true when restricted to any ``reasonable" planar regions in the first quadrant. This may be of independent interest in geometric number theory. 

We are now in position to compute $\# Y_+^\mathcal{S}(N)$. 

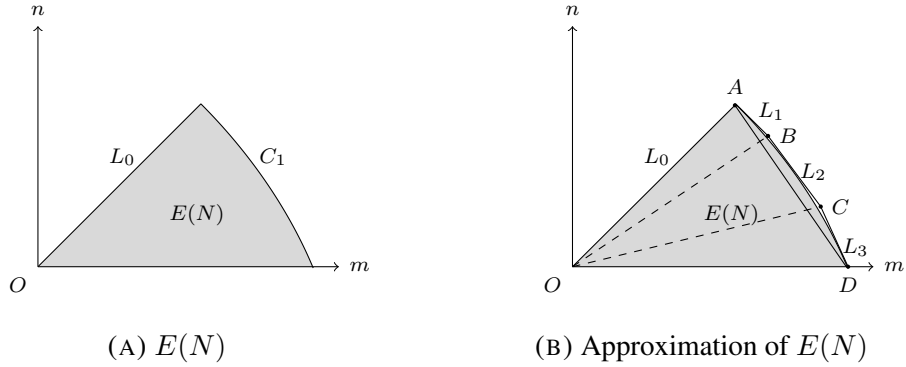
\begin{figure}[htbp]
	\centering
	
	\begin{subfigure}[b]{0.32\textwidth}
		\centering
		\begin{tikzpicture}[scale=4]
			\draw[->, thin] (0,0) -- (0,0.8) node[above] {$\scriptstyle n$};
			\draw[thin] (0,0) node[below left] {$\scriptstyle O$};		
			
			\begin{scope} 
				\clip (0,0) -- (1,1) -- (1,0) -- cycle;
				\fill[gray!30, thick, domain=0:360, smooth, variable=\t, samples=100, rotate=-45]
				plot ({1.2*cos(\t)}, {0.766*sin(\t)}) -- cycle;
			\end{scope} 
			\draw[->, thin] (0,0) -- (1.0,0) node[right] {$\scriptstyle m$};
			
			\draw[black, thin] (0,0) -- (0.542,0.542);
			\node[black, above right] at (0.2,0.3) {$\scriptstyle L_0$};
			
			\draw[black, thin, domain=57.2:90, smooth, variable=\t, samples=100, rotate=-45] 
			plot ({1.2*cos(\t)}, {0.766*sin(\t)});
			\node[black, above right] at (0.7,0.3) {$\scriptstyle C_1$};
			\node[black, above right] at (0.4,0.1) {$\scriptstyle E(N)$};
			
		\end{tikzpicture}
		\caption{$E(N)$}
		\label{dia:sub2}
	\end{subfigure}
	\qquad\qquad\qquad
	\begin{subfigure}[b]{0.32\textwidth}
		\centering
		\begin{tikzpicture}[scale=4]
			\draw[->, thin] (0,0) -- (0,0.8) node[above] {$\scriptstyle n$};
			\draw[thin] (0,0) node[below left] {$\scriptstyle O$};
			
			\begin{scope} 
				\clip (0,0) -- (1,1) -- (1,0) -- cycle;
				\fill[gray!30, thick, domain=0:360, smooth, variable=\t, samples=100, rotate=-45]
				plot ({1.2*cos(\t)}, {0.766*sin(\t)}) -- cycle;
			\end{scope} 
			\draw[black, thin, domain=57.2:90, smooth, variable=\t, samples=100, rotate=-45] 
			plot ({1.2*cos(\t)}, {0.766*sin(\t)});
			\node[black, above right] at (0.4,0.1) {$\scriptstyle E(N)$};
			\draw[->, thin] (0,0) -- (1.0,0) node[right] {$\scriptstyle m$};
			
			\draw[black, thin] (0,0) -- (0.54,0.537);
			\node[black, above right] at (0.2,0.3) {$\scriptstyle L_0$};
			
			\draw[black, thin] (0.54,0.537) -- (0.91,0);
			\draw[black, thin] (0.545,0.539) -- (0.65,0.435);
			\node[black, above right] at (0.58,0.45) {$\scriptstyle L_1$};
			\draw[black, thin] (0.825,0.2) -- (0.65,0.435);
			\node[black, above right] at (0.72,0.25) {$\scriptstyle L_2$};
			\draw[black, thin] (0.916,0) -- (0.825,0.2);
			\node[black, above right] at (0.86,0.0) {$\scriptstyle L_3$};
			\draw[black, dashed] (0,0) -- (0.65,0.435);
			\draw[black, dashed] (0,0) -- (0.825,0.2);
			\fill[black] (0.54,0.537) circle (0.2pt) node[above] {$\scriptstyle A$};
			\fill[black] (0.65,0.435) circle (0.2pt) node[right] {$\scriptstyle B$};
			\fill[black] (0.825,0.2) circle (0.2pt) node[right] {$\scriptstyle C$};
			\fill[black] (0.916,0) circle (0.2pt) node[below] {$\scriptstyle D$};
		\end{tikzpicture}
		\caption{Approximation of $E(N)$}
		\label{dia:sub3}
	\end{subfigure}
	\caption{Diagram for the elliptic region and its  approximation}\label{fig:2}
\end{figure}

\noindent\textbf{Proof of Theorem \ref{Y+}}. Recall that the set $Y_+^\mathcal{S}(N)$ consists of all coprime triples $z\in \mathcal{S}\cap D_+(N)$. Similar to that for $\# X_+^\mathcal{S}(N)$, it is sufficient to count such triples in the set $ \Delta_1\cap D_+(N) = \{(m,n) \in \Omega \ | \ \phi_1(m, n) \in D_+(N)\}$, where $\Omega$ and $\phi_1$ are defined after Lemma \ref{IT_60}.  For convenience, we denote this set by $Y_+^1(N)$. Let $E(N)$ be the region enclosed by $L_0 : m=n$, $C_1: m^2 + mn +n^2=N$, and the $m-$axis, as shown in \textsc{Figure} \ref{fig:2} (A). Due to Lemma \ref{D+},  $\# Y_+^1(N)$ is also equal to the number of coprime pairs $(m, n)$ such that $m \not\equiv n \pmod 3$ in $E(N)$. In the following, we will give upper and lower bounds for $\# Y_+^1(N)$.

Choose three straight lines which are tangent to $C_1$. They are \begin{align*} 
	L_1 &: m= -n + \frac{2\sqrt{3}}{3}\sqrt{N},    \\
	L_2 &: m= (1 - \sqrt{3})n + (\sqrt{6} - \sqrt{2})\sqrt{N},    \\
	L_3 &: m= -\frac{1}{2}n + \sqrt{N}.
\end{align*}
Let $A$, $B$, $C$, $D$ denote the intersection of $L_0$ and $L_1$, $L_1$ and $L_2$, $L_2$ and $L_3$, $L_3$ and the $m-$axis, respectively. Then \begin{flalign*}
	&A =  \left(\sqrt{\frac{N}{3}}, \sqrt{\frac{N}{3}}\right),& \\
	&B = \left(\left(\sqrt{6}-\frac{2\sqrt{3}}{3}+\sqrt{2}-2\right)\sqrt{N}, \left(-\sqrt{6}+\frac{4\sqrt{3}}{3}-\sqrt{2}+2\right)\sqrt{N}\right),\\
	&C = \left(\left(-\frac{\sqrt{6}}{3}+\frac{2\sqrt{3}}{3}-\sqrt{2}+2\right)\sqrt{N}, \left(\frac{2\sqrt{6}}{3}-\frac{4\sqrt{3}}{3}+2\sqrt{2}-2\right)\sqrt{N}\right),\\
	&D = (\sqrt{N}, 0). 
\end{flalign*}

The triangular region $\Delta_{OAD}$ is contained in $E(N)$, and $E(N)$ is contained in the pentagonal region $\Delta_{OABCD}$ enclosed by $L_0$, $L_1$, $L_2$, $L_3$ and the $m-$axis, as is shown in \textsc{Figure} \ref{fig:2} (B). Thus, $\# Y_+^1(N)$ lies between the number of coprime pairs $(m, n)$ satisfying $m \not\equiv n \pmod 3$ in $\Delta_{OAD}$ and that in  $\Delta_{OABCD}$. By Lemma \ref{TRI}, Corollary \ref{TDelta} and the remarks following them, we need to compute the areas of the above two regions. The area of $\Delta_{OAD}$ is 
$S_{OAD}=\frac{\sqrt{3}}{6}N$. The areas of $\Delta_{OAB}$, $\Delta_{OBC}$, and $\Delta_{OCD}$ are, respectively,
\begin{align*}
	&S_{OAB} = \frac{\sqrt{6}-2\sqrt{3}+3\sqrt{2}-3}{3}N,\\
	&S_{OBC} = \frac{2\sqrt{6}-4\sqrt{3}+6\sqrt{2}-6}{3}N,\\
	&S_{OCD} = \frac{\sqrt{6}-2\sqrt{3}+3\sqrt{2}-3}{3}N.
\end{align*}
Hence the area of $\Delta_{OABCD}$ is \[S_{OABCD} = S_{OAB} + S_{OBC}+S_{OCD} = \frac{4(3-\sqrt{6})(\sqrt{2}-1)}{3}N.\]
Combining with Lemma \ref{TRI} and Corollary \ref{TDelta}, we obtain \[\frac{3\sqrt{3}}{4\pi^2}N + o(N) < \# Y_+^1(N) <\frac{6(3-\sqrt{6})(\sqrt{2}-1)}{\pi^2}N + o(N).\]

\par The estimation for $\#(\Delta_3\cap D_+(N))$ is similar, and the result is the same. In summary, we have
\[\frac{3\sqrt{3}}{\pi^2}N + o(N)< \# Y_+^\mathcal{S}(N) <\frac{24(3-\sqrt{6})(\sqrt{2}-1)}{\pi^2}N + o(N),\] which is Theorem \ref{Y+}.\zb\\

\noindent\textbf{Conclusions.} In light of Theorem \ref{X+} and Theorem \ref{Y+}, we have \[\frac{\pi^2}{12(3-\sqrt{6})(\sqrt{2}-1)} \frac{1}{\log N}(1+ o(1))< \mathbb{P}_+^\mathcal{S}(N) <\frac{2\sqrt{3}\pi^2}{9}\frac{1}{\log N}(1+ o(1)),\]
where $\frac{\pi^2}{12(3-\sqrt{6})(\sqrt{2}-1)} \approx 3.6069$, and $\frac{2\sqrt{3}\pi^2}{9}\approx 3.7988$. In Theorem \ref{propP+}, we have
\[\mathbb{P}_+(N) = \frac{3\zeta(3)}{\log N}(1 + o(1)),\] with
$3\zeta(3)\approx 3.6062$. Thus we conclude that 
\begin{align*}
 \liminf \frac{\mathbb{P}_+^\mathcal{S}(N)}{\mathbb{P}_+(N)}&\geq \frac{\pi^2}{36(3-\sqrt{6})(\sqrt{2}-1)\zeta(3)}\approx 1.0002;\\
\limsup \frac{\mathbb{P}_+^\mathcal{S}(N)}{\mathbb{P}_+(N)}&\leq \frac{2\sqrt{3}\pi^2}{27\zeta(3)}\approx 1.0534,
\end{align*}
establishing the main theorem of this paper. 
 
It is worth noting that an arbitrary algebraic surface $\mathcal{A}\subseteq\mathbb{C}^3$ usually doesn't contain any integer points, and even if it contains infinitely many integer points, the ratio of the densities $\frac{\mathbb{P}_+^\mathcal{A}(N)}{\mathbb{P}_+(N)}$ is generally less than $1$. A simple example is the surface $\mathcal{A}=\{(z_0, z_1, z_2)\in \mathbb{C}^3\mid z_2=0\}$, for which  
\[\lim_{N\to \infty}\frac{\mathbb{P}_+^\mathcal{A}(N)}{\mathbb{P}_+(N)}=\frac{2\zeta(2)}{3\zeta(3)}\approx 0.9123.\]
This phenomenon leads us to suspect that the eigensurfaces of certain groups have intrinsic arithmetic structure that enables them to capture a higher density of prime tuples. Indeed, investigating potential links between groups and the density of prime triples in their eigensurfaces is at the heart of this research. We finish this discussion by raising two questions:
\begin{ques}
$ $

{\bf 1}. Does $ \lim_{N\to \infty} \frac{\mathbb{P}_+^\mathcal{S}(N)}{\mathbb{P}_+(N)}$ exist?

{\bf 2}. Which finite groups possess eigensurfaces that contain a higher density of prime tuples?
\end{ques}
We leave further discussions to another paper.\\

\noindent\textbf{Acknowledgments.} The work in this paper is partially supported by the National Natural Science Foundation of China Grant No. 12271090, Jiangsu Provincial Scientific Research Center of Applied Mathematics Grant No. BK20233002. The third author is supported partially by the Simons Foundation Grant No. MP-TSM-00002315.

\vspace{5mm}

\end{document}